\documentclass[runningheads]{llncs}
\usepackage{graphicx}
\usepackage[symbol]{footmisc}

\usepackage{float}
\usepackage{mathtools}
\usepackage{amssymb, amsmath, latexsym}
\usepackage{nicefrac}
\usepackage{hhline}
\usepackage{makecell}
\usepackage{caption}
\usepackage{multirow}
\usepackage{wrapfig}

\usepackage{colortbl}
\definecolor{bgcolor}{rgb}{0.8,1,1}
\definecolor{bgcolor2}{rgb}{0.68,0.78,0.81}
\usepackage{threeparttable}

\usepackage{pifont}
\usepackage[colorinlistoftodos,bordercolor=orange,backgroundcolor=orange!20,linecolor=orange,textsize=scriptsize]{todonotes}
\usepackage{algorithm}\usepackage{algpseudocode}

\allowdisplaybreaks
\graphicspath{ {images/} }
\def\<#1,#2>{\langle #1,#2\rangle}

\begin{document} 
\title{Optimal Data Splitting in Distributed Optimization for Machine Learning\thanks{The research of A. Beznosikov was supported by Russian Science Foundation (project No. 23-11-00229).}}
\titlerunning{Optimal Data Splitting in Distributed Optimization for ML}
%
\author{Daniil Medyakov\inst{1} \and
Gleb Molodtsov\inst{1} \and
Aleksandr Beznosikov\inst{1,2} \and
Alexander Gasnikov\inst{1,2}
}
\authorrunning{D. Medyakov, G. Molodtsov, A. Beznosikov, A. Gasnikov}
%
\institute{Moscow Institute of Physics and Technology, Moscow, Russia, \and
Institute for Information Transmission Problems, Moscow, Russia
}
\maketitle              
\begin{abstract}

The distributed optimization problem has become increasingly relevant recently. It has a lot of advantages such as processing a large amount of data in less time compared to non-distributed methods. However, most distributed approaches suffer from a significant bottleneck -- the cost of communications. Therefore, a large amount of research has recently been directed at solving this problem. One such approach uses local data similarity. In particular, there exists an algorithm provably optimally exploiting the similarity property. But this result, as well as results from other works solve the communication bottleneck by focusing only on the fact that communication is significantly more expensive than local computing and does not take into account the various capacities of network devices and the different relationship between communication time and local computing expenses. We consider this setup and the objective of this study is to achieve an optimal ratio of distributed data between the server and local machines for any costs of communications and local computations. The running times of the network are compared between uniform and optimal distributions. The superior theoretical performance of our solutions is experimentally validated.
\end{abstract}
\section{Introduction}

\subsection{Distributed optimization}

We consider optimization problems of the following form:
\begin{equation}
    \label{eq:1}
    \underset{x\in \mathbb{R}^d}{\min} ~ f(x) = \frac{1}{n} \sum \limits_{i = 1}^{n} f _i(x),
\end{equation}
where $x\in\mathbb{R}^d$ collects parameters of a statistical model to be trained, $n$ is a number of devices/nodes, and $f_i$ is an emperical risk of the devices $i$, i.e., $f_i(x) = \frac{1}{b_i}\sum_{j = 1}^{b_i} l(x, z_i^j)$ with $z_i^1, \ldots ,z_i^{b_i}$ is a set of $b_i$ samples owned by the $i$-th device and $l(x, z_i^j)$ measures mismatch between the parameter $x$ and the label of the sample $z_i^j$. This is a direct formulation of the distributed optimization problem. Nowadays, there are some reasons to consider this.

To achieve the best results in modern machine learning optimization problems, researchers and practitioners face various challenges. Dealing with modern machine learning models remains an extremely challenging task, primarily because models are trained on increasingly large datasets. Having more data in the training sample increases the robustness and generalizability of the derived model. In this case, the data is typically processed using a network of devices, i.e., collected in a distributed manner and stored in edge nodes of the network, such as in classical clustering \cite{verbraeken2020survey} and federated  \cite{konevcny2016federated,li2020federated,kairouz2021advances} learning.

Several solution methods have been proposed to solve \eqref{eq:1}. The prototype approach involves interleaving edge devices calculations (nodes $i = 2, \ldots, n$) with communications to and from the server ($i = 1$), which maintains and updates the authoritative copy of the optimization variables, eventually producing the final solution estimate. In distributed learning of complex models, the communication overhead between devices in the network often becomes a bottleneck. Such a problem makes it necessary to develop more efficient distributed learning methods, some of which have been described in \cite{konevcny2016federated,ghosh2020communication,smith2018cocoa,gorbunov2021marina}.

\subsection{Distributed optimization under similarity}

It is trendy today in machine learning to use momentum-based methods. One of the methods for solving a distributed optimization problem is the application of Nesterov acceleration \cite{nesterov2018lectures}, which is an optimal method for smooth non-distributed deterministic optimization problems. This method can be applicable to distributed networks as follows. At each iteration we calculate the gradient locally and send the results to the server. The server average the obtained gradients and make a method step. Then the number of communications is equal to the number of iterations. In this case, we obtain optimal estimates for local computations -- $\sqrt{\kappa}, \kappa = \nicefrac{L}{\mu}$, where $L$ and $\mu$ are constants of smoothness and strong convexity of target function $f$. In case $\kappa$ is small, this approach is acceptable. However, for ill-conditioned functions with a large $\kappa$, the polynomial dependence on $\kappa$ may be unsatisfactory, due to the high cost of communications. This is often the case for many empirical risk minimization (ERM) problems where the optimal regularization parameter for test predictive performance is very small.

To further improve communication complexity, we can exploit the additional structure typically found in ERM problems, known as data similarity \cite{arjevani2015communication,shamir2014communication,matsushima2014distributed}. One can define it as the difference of function gradients, i.e., $\|\nabla f_i (x) - \nabla f_j (x)\| < \delta ~~  \forall x$. But this approach is not "natural", since if the problem is not bounded, such a $\delta$ cannot exist. Let us consider for example a quadratic problem: $\nexists~ \delta: \|(A_i - A_j)x\| \leq \delta~$  if $~ x\rightarrow \infty$. Therefore, we focus on a different setting, namely on a Hessian similarity. Specifically, for all $x$ in a suitable domain of interest and all $i \neq j; ~ i,j \in \{1,\ldots,n\}$, the difference between the Hessian matrices of local losses, denoted by $\|\nabla ^2 f_i(x) - \nabla ^2 f_j(x)\|$, is bounded by $\delta$, where $\delta > 0$ measures the degree of similarity. Under this assumption, we can estimate $\delta \sim \mathcal O(\nicefrac{1}{\sqrt{N}})$, where $N$ is the sample size per device \cite{arjevani2015communication}. This setting was investigated for the first time in \cite{shamir2014communication}. After that, lower bounds were proved for this problem in \cite{arjevani2015communication}, where communication costs are proportional to $\sqrt{\nicefrac{\delta}{\mu}}$. Then for a long time researchers tried to find methods that would reach these estimates. In particular, algorithms such as \cite{tian2022acceleration,sun2022distributed,reddi2016aide,hendrikx2020statistically,beznosikov2021distributed} was obtained. In 2022, it was possible to find the optimal method which is described in \cite{kovalev2022optimal}. 

\subsection{Various communication costs and local computations} \label{eq:1.3}

In these works, the authors made an assumption that communication costs significantly more than for local computations. Moreover, in general, works on distributed optimization, not only about the Hessian similarity, made this assumption. We look at this question from a different angle, move away from fixed big communications and make 2 \textbf{assumptions}:
\begin{itemize}
    \item [1.] The devices in the network have different capacities, i.e., they perform local computations of the same amount of data for different times.
    \item [2.] The ratio of communication costs to local computation time is a variable value that can be either $\ll 1$, $\gg1$, or even $\sim 1$.
\end{itemize}
Under such assumptions, we need a new approach to the distributed optimization problem based on the already obtained optimal algorithms. This leads us to the research question of this paper:

\begin{quote}
    \textit{ Can we find such a distribution of data among the devices in the network to reduce the actual running time of the optimal algorithm} \cite{kovalev2022optimal}\textit{~for any communication costs and local computations?}
\end{quote}

In practice, networks can run for long periods of time and as a consequence, noise can occur. In other words, communication costs and device capacities are not constant values. In that way, we make one more \textbf{assumption}:
\begin{itemize}
    \item [3.] We put communication costs and device capacities as random variables, start the network operation for a long time and measure their expectation and variance. Due to the fact that the distribution of data to devices depends on the constant communication time and device power, in reality, the optimal distribution is different on account of noise. 
\end{itemize}

Therefore, this assumption raises one more research question of measuring the variance of the program running time under the optimal data distribution.

\subsection{Contributions}

In general, we can summarize our contribution as follows:
\begin{itemize}
    \item \textbf{Generalization of the computation model.} We build a general model for computing time in networks under distributed optimization. The model is based on the optimal algorithm \cite{kovalev2022optimal} and takes into account the difference in capacities of devices, in various communication costs.
    \item \textbf{Comprehensive analysis.} We pay special attention to the particular cases and obtain results for them. We consider the case where communications are too expensive and the case of inexpensive communications (not so expensive that the communication takes longer than processing all data by just one device). Moreover, we obtain results not only taking into account the  difference in time costs, but we also consider different estimates on $\delta$.
    \item \textbf{Different techniques for obtaining a solution.}  We use different techniques: Cardano's formula, upper estimates in particular cases, finding the zero of the function using the simplest numerical methods.
    \item \textbf{Decision error due to noise.} Under the third assumption, we present the theoretical error of the program running time under communication and local capacities noise.
    \item \textbf{Experiments.} We also conduct experiments confirming that with the obtained distribution it takes less time to solve the distributed problem. Besides, we also make appropriate experiments with noise in the network.
\end{itemize}

\section{Problem Statement}

Let us stay under just the first two assumptions from Section 1.3 for now. To achieve lower communication and local gradient complexity, we can refer to Algorithm \ref{alg:1} from \cite{kovalev2022optimal}. For this purpose, the function need to be represented as a sum of a smooth convex function $f_1$ and a smooth potentially non-convex function $f - f_1$. Then the algorithm is rewritten in the following form:

\begin{algorithm}
\caption{Accelerated Extragradient}\label{alg:1}
\begin{algorithmic}
\item [1:] \textbf{Input:} $x^0 = x_f^0\in\mathbb R^d$
\item [2:] \textbf{Parameters:} $\tau\in (0, 1), \eta, \theta, \alpha > 0, K\in \{1, 2, \ldots\}$
\item [3:] \textbf{for~}$k = 0, 1, 2, \ldots, K-1$\textbf{~do}
\item [4:]    $\quad\quad x_g^k = \tau x^k + (1 - \tau)x_f^k$
\item [5:]    $\quad\quad x_f^{k+1} \approx \arg\underset{x\in\mathbb R^d}{\min} [A_{\theta}^k (x) := \langle \nabla (f - f_1)(x_g^k), x - x_g^k \rangle + \frac{1}{2\theta}\|x - x_g^k\|^2 + f_1(x)]$ 
\item [6:] $\quad\quad x^{k+1} = x^k + \eta\alpha(x_f^{k+1} - x^k) - \eta\nabla f(x_f^{k+1})$
\item [7:] \textbf{end for}
\item [8:] \textbf{Output:} $x^K$
\end{algorithmic}
\end{algorithm}

We analyze the work of this algorithm, namely, find out how many operations this algorithm performs per iteration. In line 5, when solving the $\arg\min$ subproblem, one local computation is performed on devices to compute $f_i(x_g^k)$, followed by one communication to transmit these results, and additional computations on the first device to find the solution $x_f^{k+1}$. Then in line 6 there is one local computation on all devices, and one communication. We obtain an expression for the total running time of the algorithm. Let us introduce the following notations: $\tau_i$ -- time of one local computation on the $i$-th device, $K$ -- number of iterations, $\tau_{comm}$ -- time of one communication, $k_{some}$ -- additional computation of the first/central node, $n$ -- number of nodes in the network. Taking this into account, the total running time of the algorithm can be written as:
\begin{equation*}
    T_{sum} = 2\cdot\max(\tau_1, \tau_2, \ldots, \tau_n)\cdot K + 2\cdot K\cdot\tau_{comm} + \tau_1\cdot k_{some}.
\end{equation*}

Our task is to minimize the time $T_{sum}$. In view of the statement \eqref{eq:1} and the form of functions $f_i$ let us represent the time $\tau_i$ as $\tau_i = \tau_i^{loc}\cdot b_i$, where $\tau_i^{loc}$ is capacity, i.e., the time spent by the $i$-th device to process a unit of information submitted to its input, and $b_i$ is the size of dataset submitted to the $i$-th device. All $b_i$ satisfy the following constraints: $\sum\limits_{i = 1}^{n} b_i = N$, where $N$ is the size of the whole dataset, $\delta = \frac{L}{\sqrt{b_i}}$ or $\delta = \frac{L}{b_i}$ \cite{hendrikx2020statistically}.
Finally, we obtain the following optimization problem:
\begin{equation}
    \label{eq:4}
    \underset{\sum\limits_{i = 1}^{n} b_i = N; \delta = \frac{L}{{b_1}^{\gamma}}}{\min}[ 2\cdot\max(\tau_1^{loc}\cdot b_1, \ldots, \tau_n^{loc}\cdot b_n)\cdot K + 2\cdot K\cdot\tau_{comm} + \tau_1\cdot k_{some}], ~ \gamma \in \{\frac{1}{2}, 1\}.
\end{equation}

\section{How to solve (\ref{eq:4})}

\subsection{The primary problem of minimization}
In the work \cite{kovalev2022optimal} the estimates of $K$ and $k_{some}$ are presented, namely: $2\cdot K = \mathcal O(\max\{1, \sqrt{\frac{\delta}{\mu}}\}\log\frac{1}{\varepsilon}), k_{some} = \mathcal O(\max\{1, \sqrt{\frac{L}{\delta}}, \sqrt{\frac{\delta}{\mu}}, \sqrt{\frac{L}{\mu}}\}\log\frac{1}{\varepsilon})$.

Thus, \eqref{eq:4} is reduced to:

\begin{equation}
    \begin{split}
    \label{eq:5}
        \underset{\sum\limits_{i = 1}^{n} b_i = N; \delta = \frac{L}{{b_1}^\gamma}}{\min}[&(\max(\tau_1^{loc}\cdot b_1, \ldots, \tau_n^{loc}\cdot b_n) + \tau_{comm}) \cdot \mathcal O(\max\{1, \sqrt{\frac{\delta}{\mu}}\log\frac{1}{\varepsilon}\}) 
        \\&+
        \tau_1^{loc}\cdot b_1\cdot\mathcal O(\max\{1,\sqrt{\frac{L}{\delta}}, \sqrt{\frac{\delta}{\mu}}, \sqrt{\frac{L}{\mu}}\}\log\frac{1}{\varepsilon})] ,  ~ \gamma \in \{\frac{1}{2}, 1\}. 
    \end{split}
\end{equation}

\subsection{Auxiliary problem}
Consider an auxiliary problem:
\begin{equation}
    \label{eq:6}
    \underset{\sum\limits_{i = 2}^{n} b_i = N}{\min} [\max(\tau_2^{loc}\cdot b_2, \tau_3^{loc}\cdot b_3, \ldots, \tau_n^{loc}\cdot b_n)].
\end{equation}

\begin{lemma}
    \label{l1}
    The solution of problem \eqref{eq:6} is $\overrightarrow{b} = (b_2, b_3, \ldots, b_n)^{T}$ satisfying $\tau_2^{loc}\cdot b_2 = \tau_3^{loc}\cdot b_3 = \ldots = \tau_n^{loc}\cdot b_n$.
    \begin{proof}
        
        Without loss of generality, let us assume fixed values for $\tau_2^{loc}\leq \tau_3^{loc}\leq \ldots \leq \tau_n^{loc}$. 
        Then let us arbitrarily choose $b_2\geq b_3\geq \ldots \geq b_n$. 
        This is indeed the case, otherwise we would have a situation where $\exists ~ i \neq j: ~ i, j\in \{2, \ldots, n\} : \max(\tau_i^{loc}\cdot b_i, \tau_j^{loc}\cdot b_j) > \max(\tau_i^{loc}\cdot b_j, \tau_j^{loc}\cdot b_i)$, and therefore the distribution would be suboptimal. 
        
        Our goal is to minimize the function $g(\overrightarrow{b}) = \max(\tau_2^{loc}\cdot b_2, \tau_3^{loc}\cdot b_3, \ldots, \tau_n^{loc}\cdot b_n)$. 
        Suppose that there exists a distribution such that $\exists i \in \{2, \ldots, n\}: g(\overrightarrow{b}^0) = \tau_i^{loc}\cdot b_i^0$ is the minimum, and $\forall j: j \geq 2, j \neq i \hookrightarrow \tau_i^{loc}\cdot b_i^0 > \tau_j^{loc}\cdot b_j^0$. 
        It follows that $b_i^0 > \frac{\tau_j^{loc}}{\tau_i^{loc}}b_j^0 > \frac{\tau_{j_1}^{loc}}{\tau_i^{loc}}b_{j_1}^0 > \ldots > \frac{\tau_{j_k}^{loc}}{\tau_i^{loc}}b_{j_k}^0$. 
        Then, considering $\sum\limits_{i = 2}^n b_i = N \hookrightarrow b_i^0 + \frac{\tau_j^{loc}}{\tau_i^{loc}}b_j^0 + \frac{\tau_{j_1}^{loc}}{\tau_i^{loc}}b_{j_1}^0 + \ldots + \frac{\tau_{j_k}^{loc}}{\tau_i^{loc}}b_{j_k}^0 > N$,
        we obtain 
        $b_i^0 > N(1 + \tau_i^{loc}\sum\limits_{\substack{j = 2, j \neq i}}^n \frac{1}{\tau_j^{loc}})^{-1}$.
        
        Next, let us consider $b_i = N(1 + \tau_i^{loc}\sum\limits_{\substack{j = 2, j \neq i}}^n \frac{1}{\tau_j^{loc}})^{-1}, \quad b_j = \frac {\tau_i^{loc}}{\tau_j^{loc}}\cdot b_i  ~~ \forall j \in \{2,\ldots,n\}$. 
        This distribution yields a minimum of $g(\overrightarrow{b}) = \tau_i^{loc}\cdot b_i = \tau_j^{loc}\cdot b_j ~~ \forall j \in \{2,\ldots,n\}$, and $g(\overrightarrow{b}) < g(\overrightarrow{b}^0)$. 
        This contradicts the assumption of optimality. 
        Thus, for the distribution that minimizes the function $g(\overrightarrow{b}) = \max(\tau_2^{loc}\cdot b_2, \tau_3^{loc}\cdot b_3, \ldots, \tau_n^{loc}\cdot b_n)$, it holds that $\tau_2^{loc}\cdot b_2 = \tau_3^{loc}\cdot b_3 = \ldots = \tau_n^{loc}\cdot b_n$.
    \end{proof}
\end{lemma}
Let us return to the problem \eqref{eq:5}. In addition to the minimum expression already studied in \eqref{eq:6}, there are additional terms in the problem \eqref{eq:5}.  $\delta$ in \eqref{eq:5} depends on the value of $b_1$, but do not depend on $b_i, i = \overline{2, n}$. From this and Lemma \ref{l1}, it follows that in the original problem \eqref{eq:5}, the data sharing between the 2nd, 3rd, and subsequent devices should be proportional. Thus, the problem \eqref{eq:5} is reduced to a new problem with additional constraints:

\begin{equation}
    \begin{split}
    \label{eq:7}
        \min_{\substack{\sum\limits_{i = 1}^{n} b_i = N; \delta = \frac{L}{{b_1}^\gamma};  \\
        \tau_2^{loc}\cdot b_2 = \ldots = \tau_n^{loc}\cdot b_n }} 
        [&(\max(\tau_1^{loc}\cdot b_1, \ldots , \tau_n^{loc}\cdot b_n) + \tau_{comm})\cdot \mathcal O(\max\{1, \sqrt{\frac{\delta}{\mu}}\log\frac{1}{\varepsilon}\})
        \\ &+
        \tau_1^{loc}\cdot b_1\cdot\mathcal O(\max\{1,\sqrt{\frac{L}{\delta}}, \sqrt{\frac{\delta}{\mu}}, \sqrt{\frac{L}{\mu}}\}\log\frac{1}{\varepsilon})] ,  ~ \gamma \in \{\frac{1}{2}, 1\}.
    \end{split}
\end{equation}

\subsection{Define the final minimization problem}
It follows from Lemma \ref{l1} that $b_i\cdot \tau_i^{loc} = \text{const} ~ \forall i \in \overline{2, n}$.
Therefore, $$ N - b_1 = \sum\limits_{i = 2}^{n} b_i = \sum\limits_{i = 2}^{n} \frac{\tau_2^{loc}\cdot b_2}{\tau_i^{loc}} = \tau_2^{loc}\cdot b_2 \cdot \sum\limits_{i = 2}^{n} \frac{1}{\tau_i^{loc}} \Rightarrow
b_2 = \frac{N - b_1}{\tau_2 ^{loc}}(\sum\limits_{i = 2}^{n} \frac{1}{\tau_i^{loc}})^{-1}.$$

As mentioned above, we consider the case of $\delta = \frac{L}{b_1}$ and case of $\delta = \frac{L}{\sqrt{b_1}}$.
\subsubsection{3.3.1 Case of $\delta = \frac{L}{b_1}$.}

There the following relations are fulfilled:
\begin{equation}
    \notag
    \gamma = 1,~~ \mu \leq \delta \leq L \Rightarrow 
    \\
    \notag
    \begin{cases}
      2\cdot K =  \mathcal O(\sqrt{\frac{L}{\mu b_1}}\log\frac{1}{\varepsilon})\\
      k_{some} = \mathcal O(\sqrt{\frac{L}{\mu}}\log\frac{1}{\varepsilon})
    \end{cases}.
\end{equation}

Substituting this estimates into \eqref{eq:7},  the problem takes the following form:
\begin{equation}
    \notag
    \underset{\sum\limits_{i = 1}^{n} b_i = N}{\min}[(\max\{\tau_1^{loc}\cdot b_1, \tau_2^{loc}\cdot b_2\} + \tau_{comm}) \cdot \mathcal O(\sqrt{\frac{L}{\mu {b_1}}}\log\frac{1}{\varepsilon}) + \tau_1^{loc}\cdot b_1 \cdot \mathcal O(\sqrt{\frac{L}{\mu}}\log\frac{1}{\varepsilon})].
\end{equation}

As a result, leaving the only variable $b_1$ in the function we pass to the final form of minimization problem:

\begin{equation}
    \begin{split}
    \label{eq:fm1}
        \underset{0 < b_1 \leq N}{\min}\mathcal{F}(b_1) = &(\max\{\tau_1^{loc}\cdot b_1; ~(N-b_1) \cdot (\sum\limits_{i = 2}^{n} \frac{1}{\tau_i^{loc}} )^{-1}\} + \tau_{comm}) \cdot 
        \\ & \cdot
        \mathcal O(\sqrt{\frac{L}{\mu b_1}}\log\frac{1}{\varepsilon}) + \tau_1^{loc}\cdot b_1 \cdot \mathcal O(\sqrt{\frac{L}{\mu}}\log\frac{1}{\varepsilon})]. \hspace{1cm}
    \end{split}
\end{equation}

Let us investigate the problem further. To do this, find the point at which the expressions under the maximum coincide:
\begin{equation}
    \notag
    b_1^0 \cdot (\tau_1^{loc} + (\sum\limits_{i = 2}^{n} \frac{1}{\tau_i^{loc}})^{-1}) = N (\sum\limits_{i = 2}^{n} \frac{1}{\tau_i^{loc}})^{-1} \Rightarrow b_1^0 = \frac{N (\sum\limits_{i = 2}^{n} \frac{1}{\tau_i^{loc}})^{-1}}{\tau_1^{loc} + (\sum\limits_{i = 2}^{n} \frac{1}{\tau_i^{loc}})^{-1}}.
\end{equation}

Thus, we obtained two half-intervals, on each of which we can formulate a different minimization problem:
\begin{eqnarray*}
\notag
    \begin{cases}
    (a): ~ 0 < b_1 \leq b_1^0 
    \\
    (b): ~ b_1^0 <  b_1 \leq N 
    \end{cases}
\end{eqnarray*}

\begin{eqnarray*}
\notag
    \begin{cases}
    (a): ~ \max\{\tau_1^{loc}\cdot b_1; ~(N-b_1) \cdot (\sum\limits_{i = 2}^{n} \frac{1}{\tau_i^{loc}} )^{-1}\} = 
    (N-b_1) \cdot (\sum\limits_{i = 2}^{n} \frac{1}{\tau_i^{loc}})^{-1}
    \\
    (b): ~ \max\{\tau_1^{loc}\cdot b_1; ~(N-b_1) \cdot (\sum\limits_{i = 2}^{n} \frac{1}{\tau_i^{loc}} )^{-1}\} = \tau_1^{loc}\cdot b_1
    \end{cases}
\end{eqnarray*}

We construct functions of one variable $\mathcal{F}_1(b_1), \mathcal{F}_2(b_1)$ on the corresponding half-intervals that need to be minimized according to the problem \eqref{eq:fm1}:

$$
\begin{cases}
(a): ~\mathcal{F}_1(b_1) = [N (\sum\limits_{i = 2}^{n} \frac{1}{\tau_i^{loc}})^{-1} + \tau_{comm}]\cdot 
c_1 \sqrt{\frac{L}{\mu}}\log (\frac{1}{\varepsilon})  b_1^{-\frac{1}{2}} ~ - \\
\hspace{0.8cm} - ~ c_1  \sqrt{\frac{L}{\mu}}\log (\frac{1}{\varepsilon})(\sum\limits_{i =
2}^{n} \frac{1}{\tau_i^{loc}})^{-1} b_1^{\frac{1}{2}}  + \tau_1^{loc}\cdot c_2  \sqrt{\frac{L}{\mu}}\log (\frac{1}{\varepsilon}) b_1
\\
(b): ~\mathcal{F}_2(b_1) = \tau_{comm}\cdot 
c_1 \sqrt{\frac{L}{\mu}}\log (\frac{1}{\varepsilon})  b_1^{-\frac{1}{2}} + 
c_1  \sqrt{\frac{L}{\mu}}\log (\frac{1}{\varepsilon})\tau_1^{loc} b_1^{\frac{1}{2}}~ + \\ \hspace{0.8cm} + ~ \tau_1^{loc}\cdot c_2  \sqrt{\frac{L}{\mu}}\log (\frac{1}{\varepsilon}) b_1 
\end{cases}.$$

Besides we can immediately find their derivatives for further analysis:

$$
\begin{cases}
(a): ~\mathcal{F'}_1(b_1) = -\frac{1}{2}c_1 b_1^{-\frac{3}{2}}  [N (\sum\limits_{i = 2}^{n} \frac{1}{\tau_i^{loc}})^{-1} + \tau_{comm}]\cdot 
\sqrt{\frac{L}{\mu}}\log\frac{1}{\varepsilon} ~ - \\
\hspace{0.8cm} - ~
\frac{1}{2} c_1 b_1^{-\frac{1}{2}}   \sqrt{\frac{L}{\mu}}\log (\frac{1}{\varepsilon})(\sum\limits_{i = 2}^{n} \frac{1}{\tau_i^{loc}})^{-1} +
\tau_1^{loc}\cdot c_2  \sqrt{\frac{L}{\mu}}\log (\frac{1}{\varepsilon})
\\
(b): ~\mathcal{F'}_2(b_1) = -\frac{1}{2}c_1 b_1^{-\frac{3}{2}} \tau_{comm}\cdot \sqrt{\frac{L}{\mu}}\log (\frac{1}{\varepsilon}) + \frac{1}{2} c_1 b_1^{-\frac{1}{2}}  \sqrt{\frac{L}{\mu}}\log (\frac{1}{\varepsilon})\tau_1^{loc} + ~ \\ \hspace{0.6cm} ~  + \tau_1^{loc}\cdot c_2  \sqrt{\frac{L}{\mu}}\log \frac{1}{\varepsilon}
\end{cases}.
$$

\subsubsection{3.3.2 Case of $\delta = \frac{L}{\sqrt{b_1}}$.}

Here we can proceed similarly to the previous point. First, let us present the necessary relations in this case.
\begin{equation}
    \notag
    \gamma = \frac{1}{2},~~ \mu \leq \delta \leq L \Rightarrow 
    \\
    \notag
    \begin{cases}
      2\cdot K = \mathcal O(\sqrt{\frac{L}{\mu \sqrt{b_1}}}\log\frac{1}{\varepsilon})\\
      k_{some} = \mathcal O(\sqrt{\frac{L}{\mu}}\log\frac{1}{\varepsilon})
    \end{cases}.
\end{equation}

Substituting these relations into \eqref{eq:7}, we obtain

\begin{equation}
    \notag
    \underset{\sum\limits_{i = 1}^{n} b_i = N}{\min}[(\max\{\tau_1^{loc}\cdot b_1, \tau_2^{loc}\cdot b_2\} + \tau_{comm}) \cdot \mathcal O(\sqrt{\frac{L}{\mu \sqrt{b_1}}}\log\frac{1}{\varepsilon}) + \tau_1^{loc}\cdot b_1 \cdot \mathcal O(\sqrt{\frac{L}{\mu}}\log\frac{1}{\varepsilon})].
\end{equation}

Again, getting rid of all variables except $b_1$ we write the final minimization problem in this case

\begin{equation}
    \begin{split}
    \label{eq:fm2}
        \underset{0 < b_1 \leq N}{\min}[&(\max\{\tau_1^{loc}\cdot b_1; ~(N-b_1) \cdot (\sum\limits_{i = 2}^{n} \frac{1}{\tau_i^{loc}} )^{-1}\} + \tau_{comm}) \cdot \mathcal O(\sqrt{\frac{L}{\mu \sqrt{b_1}}}\log\frac{1}{\varepsilon})
        \\ &+ ~
        \tau_1^{loc}\cdot b_1 \cdot \mathcal O(\sqrt{\frac{L}{\mu}}\log\frac{1}{\varepsilon})]. \hspace{3.7cm}
    \end{split}
\end{equation}

Similarly, we select the point $b_1^0$, it turns out to be the same as in the previous paragraph. After we can obtained two half-intervals, on each of which we can formulate a different minimization problem:
\begin{eqnarray*}
    \begin{cases}
    (a): ~ 0 < b_1 \leq b_1^0 \Rightarrow \max\{\tau_1^{loc}\cdot b_1; ~(N-b_1) \cdot (\sum\limits_{i = 2}^{n} \frac{1}{\tau_i^{loc}} )^{-1}\} = ~ \\ \hspace{0.8cm} = ~
    (N-b_1) \cdot (\sum\limits_{i = 2}^{n} \frac{1}{\tau_i^{loc}})^{-1}
    \\
    (b): ~ b_1^0 <  b_1 \leq N \Rightarrow \max\{\tau_1^{loc}\cdot b_1; ~(N-b_1) \cdot (\sum\limits_{i = 2}^{n} \frac{1}{\tau_i^{loc}} )^{-1}\} = \tau_1^{loc}\cdot b_1
    \end{cases}.
\end{eqnarray*}

We construct functions of one variable $\mathcal{F}_1(b_1), \mathcal{F}_2(b_1)$ on the corresponding half-intervals that need to be minimized according to problem \eqref{eq:fm2}:

$$
\begin{cases}
    (a): ~ \mathcal{F}_1(b_1) = [N (\sum\limits_{i = 2}^{n} \frac{1}{\tau_i^{loc}})^{-1} + \tau_{comm}]\cdot 
c_1 \sqrt{\frac{L}{\mu}}\log (\frac{1}{\varepsilon})  b_1^{-\frac{1}{4}} ~ - \\
\hspace{0.8cm} - ~ c_1  \sqrt{\frac{L}{\mu}}\log (\frac{1}{\varepsilon})(\sum\limits_{i =
2}^{n} \frac{1}{\tau_i^{loc}})^{-1} b_1^{\frac{3}{4}}  + \tau_1^{loc}\cdot c_2  \sqrt{\frac{L}{\mu}}\log (\frac{1}{\varepsilon}) \cdot b_1 \\
(b): ~\mathcal{F}_2(b_1) = \tau_{comm}\cdot 
c_1 \sqrt{\frac{L}{\mu}}\log (\frac{1}{\varepsilon})  b_1^{-\frac{1}{4}} + c_1  \sqrt{\frac{L}{\mu}}\log (\frac{1}{\varepsilon})\tau_1^{loc} b_1^{\frac{3}{4}} ~ + \\
\hspace{0.8cm} + ~
 \tau_1^{loc}\cdot c_2  \sqrt{\frac{L}{\mu}}\log (\frac{1}{\varepsilon}) \cdot b_1
\end{cases}.
$$

Besides immediately find their derivatives for further analysis:
$$
\begin{cases}
(a): ~\mathcal{F'}_1(b_1) = -\frac{1}{4}c_1 b_1^{-\frac{5}{4}}  [N (\sum\limits_{i = 2}^{n} \frac{1}{\tau_i^{loc}})^{-1} + \tau_{comm}]\cdot 
\sqrt{\frac{L}{\mu}}\log \frac{1}{\varepsilon}  ~ - \\
\hspace{0.8cm} - ~ 
\frac{3}{4} c_1 b_1^{-\frac{1}{4}}   \sqrt{\frac{L}{\mu}}\log (\frac{1}{\varepsilon})(\sum\limits_{i = 2}^{n} \frac{1}{\tau_i^{loc}})^{-1} +
\tau_1^{loc}\cdot c_2  \sqrt{\frac{L}{\mu}}\log \frac{1}{\varepsilon}
\\
(b): ~\mathcal{F'}_2(b_1) = -\frac{1}{4}c_1 b_1^{-\frac{5}{4}} \tau_{comm}\cdot \sqrt{\frac{L}{\mu}}\log \frac{1}{\varepsilon} + \frac{3}{4} c_1 b_1^{-\frac{1}{4}}  \sqrt{\frac{L}{\mu}}\log (\frac{1}{\varepsilon})\tau_1^{loc}   ~ + \\
\hspace{0.8cm} + ~ \tau_1^{loc}\cdot c_2  \sqrt{\frac{L}{\mu}}\log \frac{1}{\varepsilon}
\end{cases}.
$$

\subsection{Final solution}

\subsubsection{3.4.1 Case of $~\delta = \frac{L}{b_1}$.} \label{eq:3.4.1}

Our goal is to find the minimum of the already obtained functions $\mathcal{F}_1(b_1), \mathcal{F}_2(b_1)$. To do this, we will look for the zeros of $\mathcal{F'}_1(b_1), \mathcal{F'}_2(b_1)$. Here we obtain the cubic equation.
To solve it, we can use the Cardano's formula. Consider the equation $ax^{-\frac{1}{2}} + bx^{-\frac{3}{2}} + c = 0$,
where in cases $(a): ~ 0 < b_1 \leq b_1^0 $ and $(b): ~ b_1^0 <  b_1 \leq N$ we put:

\begin{eqnarray}
\notag
    \begin{cases}
    (a): ~ ~ a = \frac{1}{2} c_1 \sqrt{\frac{L}{\mu}}\log (\frac{1}{\varepsilon})(\sum\limits_{i = 2}^{n} \frac{1}{\tau_i^{loc}})^{-1}; \\ \hspace{1cm}
b = -\frac{1}{2} c_1 [N (\sum\limits_{i = 2}^{n} \frac{1}{\tau_i^{loc}})^{-1} + \tau_{comm}]\cdot 
\sqrt{\frac{L}{\mu}}\log \frac{1}{\varepsilon};~ c = \tau_1^{loc}\cdot c_2  \sqrt{\frac{L}{\mu}}\log (\frac{1} {\varepsilon})
    \\
\notag
    (b): ~ ~ a = \frac{1}{2} c_1  \sqrt{\frac{L}{\mu}}\log (\frac{1}{\varepsilon})\tau_1^{loc}; \\ \hspace{1cm}
b = -\frac{1}{2}c_1 \tau_{comm}\cdot \sqrt{\frac{L}{\mu}}\log \frac{1}{\varepsilon}; ~
c = \tau_1^{loc}\cdot c_2  \sqrt{\frac{L}{\mu}}\log (\frac{1} {\varepsilon})
    \end{cases}
\end{eqnarray}

Then on the condition that 

\begin{align*}
    N \geq &\frac{a^2}{3 c^2}+\frac{\sqrt[3]{2 a^6+3 \sqrt{3} \sqrt{4 a^3 b^3 c^6+27 b^4 c^8}+18 a^3 b c^2+27 b^2 c^4}}{3 \sqrt[3]{2} c^2} \\
    &- \frac{\sqrt[3]{2}\left(-a^4-6 a b c^2\right) }{3 c^2 \sqrt[3]{2 a^6+3 \sqrt{3} \sqrt{4 a^3 b^3 c^6+27 b^4 c^8}+18 a^3 b c^2+27 b^2 c^4}},
\end{align*}

we get a solution:

\begin{align*}
     x=&\frac{a^2}{3 c^2}+\frac{\sqrt[3]{2 a^6+3 \sqrt{3} \sqrt{4 a^3 b^3 c^6+27 b^4 c^8}+18 a^3 b c^2+27 b^2 c^4}}{3 \sqrt[3]{2} c^2}\\ 
     &-\frac{\sqrt[3]{2}\left(-a^4-6 a b c^2\right)} 
    {3 c^2 \sqrt[3]{2 a^6+3 \sqrt{3} \sqrt{4 a^3 b^3 c^6+27 b^4 c^8}+18 a^3 b c^2+27 b^2 c^4}}.
\end{align*}

Hence the desired solution is trivially obtained. Since we have obtained one value of $b_1$ on each of the half-intervals, which is the minimum of the function on its, so by choosing the one on which the function is smaller, we obtain the optimal value of $b_1$.

\subsubsection{3.4.2 Case of $~\delta = \frac{L}{\sqrt{b_1}}$.}\label{eq:3.4.2}

Proceed similarly as in the previous paragraph does not work, since we cannot write out the solution of these equations in analytic form due to their powers.
Therefore, let us consider the following particular cases:
\begin{enumerate}
    \item $\forall i\hookrightarrow \tau_{comm} \ll \tau_i^{loc};$
    \item $\forall i\hookrightarrow \tau_{comm} \gg \tau_i^{loc}, \forall i\neq j\hookrightarrow \tau_i^{loc} = \tau_j^{loc}.$
\end{enumerate}

Let us introduce new notation: $\alpha = c_1\cdot\sqrt{\frac{L}{\mu}}\cdot \log\frac{1}{\varepsilon},\beta = c_2\cdot\sqrt{\frac{L}{\mu}}\cdot \log\frac{1}{\varepsilon} $, for simplicity. Now we are ready to consider two cases separately. 

\textbf{Case 1:}
    \begin{itemize}
    \item [(a):] $0 < b_1 \leq b_1^0$ and 
    $~\mathcal{F}_1(b_1) = [N (\sum\limits_{i = 2}^{n} \frac{1}{\tau_i^{loc}})^{-1} + \tau_{comm}]\cdot 
    \alpha  b_1^{-\frac{1}{4}} ~ - \\ -
    \alpha(\sum\limits_{i =
    2}^{n} \frac{1}{\tau_i^{loc}})^{-1} b_1^{\frac{3}{4}}  + \tau_1^{loc}\cdot\beta b_1$. Let us assume that
    $
    \tau_1^{loc} \leq \tau_2^{loc} \leq\ldots \leq \tau_n^{loc}.
    $
    Using this assumption and $\tau_{comm} \ll \tau_i^{loc}$, one can obtain the following estimate:

    \begin{eqnarray}
        \label{eq:10}
            (\sum\limits_{i = 2}^n \frac{1}{\tau_i^{loc}})^{-1} &=& \frac{1}{\frac{1}{\tau_1^{loc}} + \ldots + \frac{1}{\tau_n^{loc}}}
            \notag\\ &=& 
            \frac{\tau_2^{loc}\cdot \ldots \cdot\tau_n^{loc}}{\tau_3^{loc}\cdot \ldots \cdot\tau_n^{loc} + \tau_2^{loc}\cdot \tau_4^{loc}\cdot\ldots \cdot\tau_n^{loc} + \ldots + \tau_2^{loc}\cdot \ldots \cdot\tau_{n-1}^{loc}} 
            \notag\\ &\geq&\frac{\tau_2^{loc}}{n - 1} \gg \tau_{comm}.
    \end{eqnarray}
    
    Given the estimate \eqref{eq:10}, the functions $\mathcal{F}_1(b_1)$ and accordingly $\mathcal{F'}_1(b_1)'$ can be approximately simplified as follows:
    \begin{eqnarray*}
            \mathcal{F}_1(b_1) &=& \alpha(\sum\limits_{i = 2}^n \frac{1}{\tau_i^{loc}})^{-1}\cdot b_1^{-\frac{1}{4}}(N - b_1) + \tau_1^{loc}\beta\cdot b_1,
            \\
        \mathcal{F'}_1 (b_1) &=& \alpha(\sum\limits_{i = 2}^n \frac{1}{\tau_i^{loc}})^{-1}\cdot (-\frac{1}{4}b_1^{-\frac{5}{4}}N - \frac{3}{4}b_1^{-\frac{1}{4}}) + \tau_1^{loc}\beta.
    \end{eqnarray*}
    
    We get the equation in the same powers, and then again we cannot write out an analytic solution, but for this problem it is easier to find a numerical solution.
    
    \item[(b):] $b_1^0\leq b_1\leq N$ and $\mathcal{F}_2(b_1) = \tau_{comm}\cdot 
    \alpha  b_1^{-\frac{1}{4}} + 
    \alpha\tau_1^{loc} b_1^{\frac{3}{4}}  + \tau_1^{loc}\cdot \beta b_1 ~ = \\ = ~ \alpha\cdot b_1^{-\frac{1}{4}}(\tau_{comm} + \tau_1^{loc}b_1) + \beta\cdot\tau_1^{loc}\cdot b_1$. With the same assumption that $
    \tau_1^{loc} \leq \tau_2^{loc} \leq\ldots \leq \tau_n^{loc}
    $, we get
    
    \begin{equation}
    \label{eq:11}
    \begin{split}
      \tau_1^{loc}b_1 {\geq} \frac{\tau_1^{loc}N\frac{\tau_2^{loc}}{n - 1}}{\tau_1^{loc} + \frac{\tau_n^{loc}}{n - 1}}\geq \frac{\tau_1^{loc}\tau_2^{loc}N}{(n - 1)(\tau_1^{loc} + \tau_n^{loc})}\geq \frac{\tau_1^{loc}\tau_2^{loc}N}{2(n - 1)\tau_n^{loc}} ~ \gg \\
       \gg ~ \tau_{comm}\frac{N}{2(n - 1)} \gg \tau_{comm} \hspace{3cm}
    \end{split}
    \end{equation}
    
    Here, using \eqref{eq:11}, we can also simplify $F_2(b_1)$ and then $\mathcal{F'}_2(b_1)$:
    \begin{eqnarray*}
        \mathcal F_2(b_1) = \alpha\cdot\tau_1^{loc}\cdot b_1^{\frac{3}{4}} + \beta \tau_1^{loc}\cdot b_1, \quad 
        \mathcal{F'}_2(b_1) = \frac{3}{4}\alpha\cdot\tau_1^{loc\cdot} b_1^{-\frac{1}{4}} + \beta\cdot\tau_1^{loc} > 0.
    \end{eqnarray*}
    
    Since the derivative of the function is positive, the function is increasing, and therefore the minimum is taken at $b_1 = b_1^{0} = \frac{N (\sum\limits_{i = 2}^{n} \frac{1}{\tau_i^{loc}})^{-1}}{\tau_1^{loc} + (\sum\limits_{i = 2}^{n} \frac{1}{\tau_i^{loc}})^{-1}}.$  
Thus, in the case of small $\tau_{comm}$ we obtained the following result: 

\begin{equation}
    \label{eq:temp505}
    b_{1,{\min}} \leq b_1^0 = \frac{N (\sum\limits_{i = 2}^{n} \frac{1}{\tau_i^{loc}})^{-1}}{\tau_1^{loc} + (\sum\limits_{i = 2}^{n} \frac{1}{\tau_i^{loc}})^{-1}}.
\end{equation}
\end{itemize}

\textbf{Case 2:}
\begin{itemize}
\item []
Here we also define: $\tau := \tau_i^{loc}~ \forall i \in {1,\ldots, n} $. Then we can rewrite the target function of \eqref{eq:fm1} in the following way:
\begin{equation*}
        \mathcal{F}(b_1) = (\max\{\tau b_1; (N-b_1) \frac{\tau}{n-1}\} + \tau_{comm}) \cdot \frac{\alpha}{\sqrt[4]{b_1}}+\tau \beta b_1. 
\end{equation*}

Consider the case $\tau_{comm} = N^2 \tau$. We can assume also that the size of data $N$ is large, therefore $\tau_{comm} \gg N\tau$. And then:
\begin{equation}
    \notag
     \max \{\tau b_1; (N-b_1)\frac{\tau}{n-1}\} < \tau N \ll \tau_{comm}\Rightarrow \mathcal{F}(b_1) \approx \frac{\alpha \tau_{comm}}{\sqrt[4]{b_1}} + \beta \tau b_1 
\end{equation}
\begin{equation}
    \notag
    \mathcal{F}' (b_1) = -\frac{\alpha \tau_{comm}}{4b_1\sqrt[4]{b_1}} + \beta \tau = 0 \Rightarrow b_{1, {\min}}^\frac{5}{4} = \frac{\tau _{comm}\alpha}{4\beta\tau}\Rightarrow b_{1, {\min}} = (\frac{\tau _{comm}\alpha}{4\beta\tau})^{\frac{4}{5}}.
\end{equation}
If the found value $b_{1, {\min}}$ lies in the interval $(0, N) $, one can found the optimal value of $\mathcal{F}$:
\begin{equation}
    \notag
    \begin{split}
    \mathcal{F}(b_{1_{\min}}) = (\alpha \tau_{comm})^\frac{4}{5} \cdot (4\beta\tau)^\frac{1}{5} + (\beta \tau)^\frac{1}{5}(\frac{\alpha \tau_{comm}}{4})^\frac{4}{5} ~ = \\
    = ~ (\alpha \tau_{comm})^\frac{4}{5}(\beta\tau)^\frac{1}{5}(4^\frac{1}{5} + 4^{-\frac{4}{5}}). \hspace{2cm}
    \end{split}
\end{equation}
Otherwise, the minimum is reached at the right boundary, since at zero we can say that the function increases. Summarizing all of the above in this case, it is worth noting that for large values of $N$ the second special case generalizes to the following condition:
\begin{equation}
    \label{eq:temp404}
    \begin{split}
    \forall i  \hookrightarrow \tau_{comm} = \mathcal{O}( N^k \tau_i^{loc}) ~\text{with}~ k >1 , ~\text{and}~ \forall i\neq j\hookrightarrow \tau_i^{loc} = \tau_j^{loc} = \tau,
    \\
    \min {\mathcal{F}}(b_1) = \begin{cases}
      (\alpha \tau_{comm})^\frac{4}{5}\cdot (\beta\tau)^\frac{1}{5}(4^\frac{1}{5} + 4^{-\frac{4}{5}}), & 0 < (\frac{\tau _{comm}\alpha}{4\beta\tau})^{\frac{4}{5}} < N\\
      \frac{\alpha\tau _{comm}}{N} + \beta \tau N , & (\frac{\tau _{comm}\alpha}{4\beta\tau})^{\frac{4}{5}} \geq N
    \end{cases}.
    \end{split}
\end{equation}
\end{itemize}

\subsection{Numerical solution}\label{eq:3.5}

Since an analytical solution is not found for all cases, we can give a general numerical solution to our problem. In order to determine the minimum of these functions on the respective half-intervals, we examine points where the derivatives of $\mathcal{F'}_1(b_1)$ and $\mathcal{F'}_2(b_1)$ approach zero. It should be noted that, given the nature of these functions, their derivatives can only be zero once on the desired half-interval. Hence, by employing the Newton's method \cite{polyak2007newton} for $\mathcal{F'}_1(b_1)$ and $\mathcal{F'}_2(b_1)$, we can locate its zeros. Subsequently, we need to compare the values of the corresponding function at these points with the value at the extreme point of the interval. One of these points provides the optimal value, thereby serving as the ultimate solution to the problem \eqref{eq:fm1} and \eqref{eq:fm2}.

\section{Noise in the networks}

Now we proceed to the third assumption from Section 1.3. As mentioned above, let $\tau_{comm}$ and $\tau_i^{loc}$ be random variables with $\mathbb E [\tau_{comm}] < \infty$, $ \mathbb E [\tau_i^{loc}] < \infty$, $\mathbb D [\tau_{comm}] < \infty$, $\mathbb D [\tau_i^{loc}] < \infty$. Here we consider the case of $\delta = \frac{L}{\sqrt{b_1}}$. The analytical solution in this setup was obtained in the particular cases for small and large communication times. Let us consider them separately. 

\subsection{Case of big communication time}\label{s:4.1}

Let us consider the case of \eqref{eq:temp404}. In particular, we obtained that
\begin{equation}
    \label{eq:17}
    \mathcal{F}(b_{1, {\min}}) = (\alpha\cdot\tau_{comm})^{\nicefrac{4}{5}}\cdot (\beta\cdot\tau_1^{loc})^{\nicefrac{1}{5}}\cdot(4^{\nicefrac{1}{5}} + 4^{-\nicefrac{4}{5}})
\end{equation}
with $\alpha = c_1\cdot\sqrt{\frac{L}{\mu}}\cdot \log\frac{1}{\varepsilon},\beta = c_2\cdot\sqrt{\frac{L}{\mu}}\cdot \log\frac{1}{\varepsilon} $.

To estimate the variance of the result $\mathcal{F}(b_{1_{\min}}) $, we produce the following equality: $X, Y -$  independent random variables $\Rightarrow \mathbb D [XY] = \mathbb E [(XY - \mathbb E[XY])^2] = \mathbb E[(XY)^2] - 2\mathbb E^2 [XY] + \mathbb E^2[XY] = \mathbb E[X^2]\mathbb E[Y^2] - \mathbb E^2 [X]\mathbb E^2 [Y] = (\mathbb D [X] + \mathbb E^2 [X])\cdot(\mathbb D [Y] + \mathbb E^2 [Y]) - \mathbb E^2 [X]\mathbb E^2 [Y] = \mathbb D[X]\mathbb D[Y] + \mathbb D[X]\mathbb E^2 [Y]  + \mathbb D[Y]\mathbb E^2 [X]$. Applying this equality to \eqref{eq:17}, we get the variance of $\mathcal{F}(b_{1_{\min}}) $:

\begin{equation*}
    \begin{split}
        \mathbb D [\mathcal{F}(b_{1,{\min}})] =& [\alpha^{\nicefrac{4}{5}}\cdot\beta^{\nicefrac{1}{5}}\cdot(4^{\nicefrac{1}{5}} + 4^{-\nicefrac{4}{5}})]\cdot\{\mathbb D [(\tau_{comm})^{\nicefrac{4}{5}}]\mathbb D [(\tau_1^{loc})^{\nicefrac{1}{5}}]
        \\ &+
         \mathbb D [(\tau_{comm})^{\nicefrac{4}{5}}]\mathbb E^2 [(\tau_1^{loc})^{\nicefrac{1}{5}}] + \mathbb D [(\tau_1^{loc})^{\nicefrac{1}{5}}]\mathbb E^2 [(\tau_{comm})^{\nicefrac{4}{5}}]\}.  \hspace{1.5cm}
    \end{split}
\end{equation*}

\subsection{Case of small communication time}

Here we consider noise only in communications, i.e. we put the time of communications as a random variable with mathematical expectation and finite variance, and the time of local computations for each device as a constant value.  We consider the function value \eqref{eq:fm2} at the point \eqref{eq:temp505}:

\begin{equation*}
    \begin{split}
        \mathcal{F}(b_{1, {\min}}) = [(N - b_1^0)\cdot(\sum\limits_{i = 2}^n \frac{1}{\tau_i^{loc}})^{-1} + \tau_{comm}]\cdot\alpha\cdot\frac{1}{(b_1^0)^{\nicefrac{1}{4}}} + \tau_1^{loc}\cdot b_1^0\cdot\beta ~ = \hspace{0.4cm}
        \\ \quad =
        \left[\frac{\tau_1^{loc}\cdot N\cdot (\sum\limits_{i = 2}^n \frac{1}{\tau_i^{loc}})^{-1}}{\tau_1^{loc} + (\sum\limits_{i = 2}^n \frac{1}{\tau_i^{loc}})^{-1}} + \tau_{comm}\right]\cdot\alpha\cdot\frac{(\tau_1^{loc} + (\sum\limits_{i = 2}^n \frac{1}{\tau_i^{loc}})^{-1})^{\nicefrac{1}{4}}}{(N\cdot (\sum\limits_{i = 2}^n \frac{1}{\tau_i^{loc}})^{-1})^{\nicefrac{1}{4}}} + \tau_1^{loc}\cdot b_1^0\cdot\beta
    \end{split}
\end{equation*}
with $\alpha = c_1\cdot\sqrt{\frac{L}{\mu}}\cdot \log\frac{1}{\varepsilon},\beta = c_2\cdot\sqrt{\frac{L}{\mu}}\cdot \log\frac{1}{\varepsilon} $.

Then the required variance is as follows:
\begin{equation*}
    \label{eq:20}
    \mathbb D [\mathcal{F}(b_{1, {\min}})] = \left[\alpha\cdot\frac{(\tau_1^{loc} + (\sum\limits_{i = 2}^n \frac{1}{\tau_i^{loc}})^{-1})^{\nicefrac{1}{4}}}{(N\cdot (\sum\limits_{i = 2}^n \frac{1}{\tau_i^{loc}})^{-1})^{\nicefrac{1}{4}}}\right]^2\cdot \mathbb D [\tau_{comm}].
\end{equation*}

\section{Experiments}

\subsection{Experiments with data distribution}
For experimental verification of the theoretical results we consider the ridge regression problem: 
\begin{equation}
    \label{ridge}
    \underset{\omega}{\min}[ \frac{1}{2N} \|X\omega - y\|_2^2 + \frac{\lambda}{2}\|\omega\|_2^2],
\end{equation}
where $\omega$ is the vector of weights of the model, $\{x_i, y_i\}_{i = 1}^N$ is the training dataset, and $\lambda > 0$ is the regularization parameter. We consider a network  with 21 workers simulated on a single-CPU machine. We use dataset from LIBSVM library \cite{chang2011libsvm}. Value $\tau_1^{loc} = 1$, values for other i $\tau_i^{loc}, i \neq 1$ were taken conditionally and generated uniformly from 3 to 7. $\tau_{comm}$ were chosen so that $\frac{\tau _{comm}}{\tau _1^{loc}} = 10^l, l = \overline{-6, 12}$

\begin{remark}
In the first stages $\arg\underset{x}{\min} [p(x_k^g) + <\nabla p(x_k^g), x - x_k^g> + \frac{1}{2\Theta}||x - x_k^g||^2 + q(x)]$ was searched explicitly (Line 5 of Algorithm 1). This was done by equating the gradient to zero. Applying to the problem \ref{ridge}: 
\begin{center}
$\nabla p(\omega_k^g) + \frac{1}{\Theta}(x - \omega_k^g) + \nabla q(x) = 0 \Rightarrow \lambda \omega_k^g + \frac{1}{\Theta}(Ix - \omega_k^g) + \frac{1}{N}X^T(Xx - y) = 0 \Rightarrow$

$ x = (I \frac{1}{\Theta} + \frac{1}{N}X^TX)^{-1}(\frac{1}{\Theta} \omega_k^g + \frac{1}{N}X^Ty - \lambda \omega_k^g) $
\end{center}

Then in Algorithm 1 on line 5: $x_f^{k+1} = x$
\end{remark}
 We implement Algorithm \ref{alg:1} in Python 3.9.6 using the iterative OGM-G method from \cite{kim2021optimizing} to find the $\arg\min$ in \ref{alg:1} (that is what the original article \cite{kovalev2022optimal} recommends). After calculating the required number of iterations to achieve a certain accuracy we find the values of constants $c_1, c_2$, and, respectively, $\alpha, \beta$. With their help, we are able to distribute the data from the dataset to the devices according to the above formulas. 

Next, we run the algorithm and measure the running time on the resulting distribution of data across devices and uniform distribution. Our goal is find the acceleration between our choice of data distribution and uniform splitting.

\begin{figure}[!ht]
    {\includegraphics[scale = 0.18]{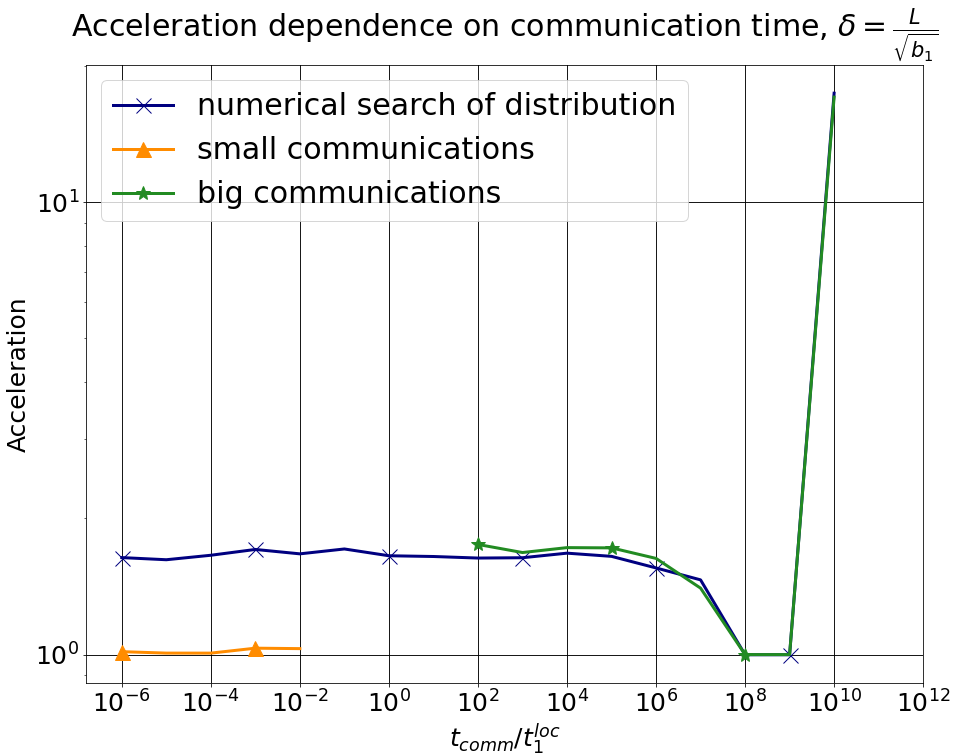}}
    {\includegraphics[scale = 0.18]{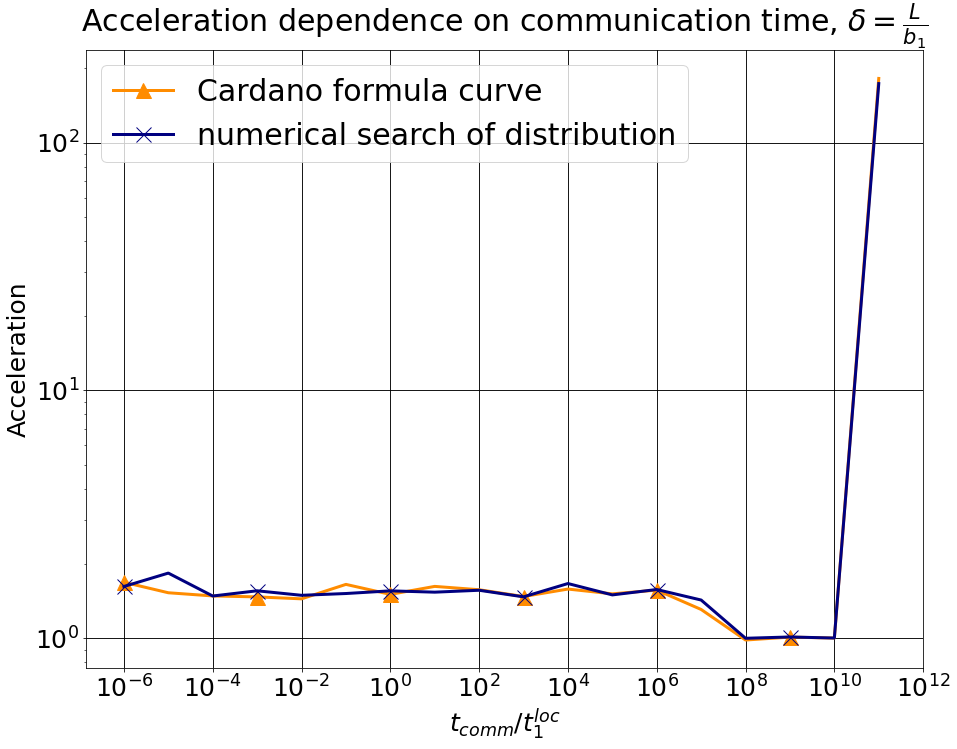}}
    \caption{Experiments with data distribution}
    \label{ris:image}
\end{figure}

Two cases of different $\delta$: $\delta = \frac{L}{\sqrt{b_1}}$, $\delta = \frac{L}{b_1}$, are considered. For the case with $\delta = \frac{L}{\sqrt{b_1}}$ we use following approaches to find $b_{1, \min}$: 1) for all cases of the communication time we use the Newton's method to find the solution numerically; 2) for small and large communications we also use results of Section 3.4.2. 
For the case with $\delta = \frac{L}{b_1}$ to find $b_{1, \min}$ we also use the Newton's method and additionally the Cardano's formula from Section 3.4.1. See Figure \ref{ris:image} for results.


Let us analyze the obtained plots. The formula for the case of large communications (Section 3.4.2. ,Case 2) and the Cardano's formula (Section 3.4.1.) practically coincided with the optimal solution search by the Newton's method. The case of small communications showed worse results. This is explained by the fact that the formula was obtained in rough approximation. 



\subsection{Experiments with noise}

\begin{figure}[h]
    \center{{\includegraphics[scale = 0.33]{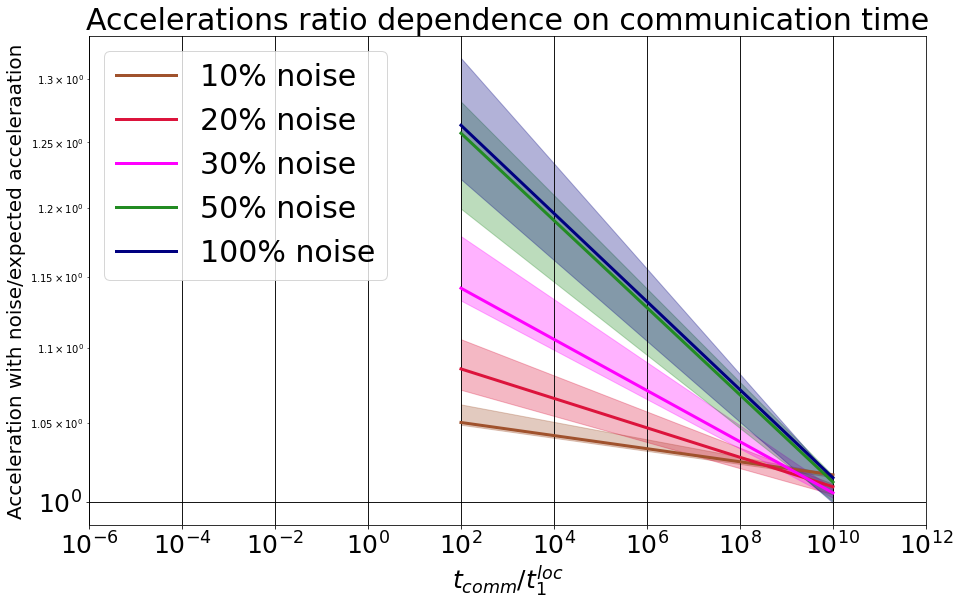}}}
    \caption{Experiments with the noise in the network}
    \label{ris:image2}
\end{figure}

We modify the simulation of Algorithm \ref{alg:1} by adding noise to communication and device power times. We generate the noise from a uniform distribution and its values are 10, 20, 30, 50, and 100 percent, respectively, relative to the absolute value of communication time and device power. In the new noise model, we measure the running time of the ridge regression problem with the resulting data distribution and with the uniform one, obtaining the acceleration that gives our data distribution. In the process, we measure the mathematical expectation of communication costs and device powers and obtained the acceleration at these expected values. Experiments are conducted only in the case of big communication costs. In Figure \ref{ris:image2}, we plots the ratio of these accelerations and also confidence intervals.

Let us analyze the obtained results. The figure shows that all the lines fall within the confidence intervals, which means that for any value of noise, the theoretical calculations given in Section 4.1 are confirmed by experiments. We also note that near the value $t_{comm}/t_1^{loc} = 10^{10}$, noise practically ceases to affect the results.

\section{Conclusion}

In this article, we presented a new data partioning method for a distributed optimization problem with different time costs of local computations and communications.e. Our solution is based on constructing the running time function of Algorithm \ref{alg:1} and finding its minimum. Our method works well in networks with varying communication costs between the server and local devices and different capacities of the devices. The theoretical results confirmed experimentally. This shows that our method gives acceleration on this type of problems. In addition, by assuming noise in the networks, we obtained the error of the optimal solution and conducted appropriate experiments.

\bibliographystyle{splncs04}
\bibliography{refs}



\end{document}